\documentclass[thmsa,11pt]{article}
\usepackage[utf8]{inputenc}
\usepackage{amsmath}
\usepackage{amssymb}

\begin{document}

\title{Corrigendum \& Addendum to:\\
\textit{Variations on a Visserian theme}}
\author{Ali Enayat}
\maketitle

\section{Introduction}

This note complements my paper \cite{Enayat-Visserian}; it was prompted by
the discovery of a gap by Piotr Gruza and Mateusz \L e\l yk in the proof of
solidity of $\mathsf{KM}$ (Kelley-Morse theory of classes), as presented in 
\cite[Theorem 2.7]{Enayat-Visserian}. At the time of this writing it is open
whether $\mathsf{KM}$ and its higher order variants are solid theories, but
as explained in the proof of Theorem 2.7,, the proof presented in \cite%
{Enayat-Visserian} works fine for the strengthening of $\mathsf{KM}$ with
the fragment $\mathsf{CC}_{\mathsf{set}}$ of the class choice scheme\textit{%
\ }(see Definition 2.2).\footnote{%
As noted in the first paragraph of the proof of Theorem 2.7, there is no
need for $\mathsf{CC}$ if the interpretations at work are identity
preserving.}\textit{\ }The solidity of higher order variants of $\mathsf{KM}%
_{n}$ augmented with appropriate analogues of $\mathsf{CC}_{\mathsf{set}}$
can be similarly established, as indicated in Remark 2.6. We assume that the
reader has \cite{Enayat-Visserian} for ready reference. \medskip 

It is also worth noting that there has been a number of advances in relation
to the topics investigated in \cite{Enayat-Visserian}; see \cite%
{Enayat-Lelyk-JPM} and \cite{Enayat-Lelyk-correction} for the relevant
references.\medskip

\textbf{Acknowledgements.} I am grateful to Piotr Gruza and Mateusz \L e\l %
yk for pointing out the gap in the proof of \cite[Theorem 2.7]%
{Enayat-Visserian}. Hats off also to Vika Gitman for her comments on an
earlier draft, which have been instrumental in shaping the proof of Theorem
2.7.

\section{Solidity of KM and its higher order variants}

The proof of solidity of $\mathsf{KM}$, as presented in \cite[Theorem 2.7]%
{Enayat-Visserian} has two stages. In the first, one repeats the proof
strategy of the proof of solidity of $\mathsf{ZF}$ to build a certain
isomorphism between the first sorts of two models of $\mathsf{KM}$; and in
the second stage the isomorphism built in the first stage is naturally
lifted to an isomorphism to the objects in the second sort. As we will see,
there is a subtle gap in the first stage of the proof of \cite[Theorem 2.7]%
{Enayat-Visserian} that can be circumvented if $\mathsf{KM}$ is strengthened
to $\mathsf{KM+CC}_{\mathsf{set}}$. The second stage is fine as is, and in
particular there is no need for $\mathsf{CC}_{\mathsf{set}}$ there. \medskip

We take this occasion to describe a more perspicuous variant of the proof of
solidity of $\mathsf{ZF}$ presented in \cite{Enayat-Visserian}. We will then
use this new presentation to pinpoint the role of $\mathsf{CC}_{\mathsf{set}}
$ in fixing the gap in the proof of solidity of $\mathsf{KM}$ in \cite[%
Theorem 2.7]{Enayat-Visserian}. \medskip 

\noindent \textbf{2.1.~Theorem.} $\mathsf{ZF}$ \textit{is solid}.\medskip

\noindent \textbf{Proof.} Suppose $\mathcal{M},\mathcal{\ M}^{\ast }$, and $%
\mathcal{N}$ are models of $\mathsf{ZF}$ such that $\mathcal{M}%
\trianglerighteq _{\mathrm{par}}\mathcal{N}\trianglerighteq _{\mathrm{par}}%
\mathcal{M}^{\ast }$, and there is an $\mathcal{M}$-definable isomorphism%
\textit{\ }$i_{0}:\mathcal{M}\rightarrow \mathcal{M}^{\ast }.$ Using the
`Scott-trick', we can assume without loss of generality, that the above
interpretation of $\mathcal{N}$ in $\mathcal{M}$, and $\mathcal{M}^{\ast }$
in $\mathcal{N}$, are both \textit{identity preserving}.\footnote{%
In model-theoretic language, this feature of $\mathsf{ZF}$ is phrased as: $%
\mathsf{ZF}$\textit{\ eliminates imaginaries}. However, many foundational
theories, including $\mathsf{KM}$, fail to eliminate imaginaries; see
Theorem 2.4.} We will show that there is an $\mathcal{N}$-definable
isomorphism between $\mathcal{N}$ and $\mathcal{M}^{\ast }$. Since $\mathcal{%
M}$ injects $M$ into $M^{\ast }$ via $i_{0}$, and $M^{\ast }\subseteq N$, we
have:$\medskip $

\noindent (1) \ \ $N$ is a proper class as viewed from $\mathcal{M}$.
\medskip

\noindent Thanks to the assumptions that $\mathcal{M}\trianglerighteq _{%
\mathrm{par}}\mathcal{N}$ and $\mathcal{M}\models \left[ i_{0}:(\mathbf{V}%
,\in )\overset{\cong }{\longrightarrow }\mathcal{M}^{\ast }\right] $, it is
easy to see that:$\medskip $

\noindent (2) $\ \ \mathcal{N}$ views $\in ^{\mathcal{M}^{\ast }}$as
well-founded in the strong sense that if $X$ is a non-empty subset of $%
M^{\ast }$ that is $\mathcal{N}$-definable, then $X$ has an $\in ^{\mathcal{M%
}^{\ast }}$-minimal element.$\medskip $

\noindent Next, let

\begin{equation*}
I=\left\{ \alpha \in \mathbf{Ord}^{\mathcal{N}}\mid \mathcal{N}\models \left[
\exists \beta \in \mathbf{Ord}^{\mathcal{M}^{\ast }}\ \exists f:(V_{\alpha
},\in )\overset{\cong }{\longrightarrow }(V_{\beta },\in )^{\mathcal{M}%
^{\ast }}\right] \right\} ,
\end{equation*}%
and 
\begin{equation*}
J=\left\{ \beta \in \mathbf{Ord}^{\mathcal{M}^{\ast }}\mid \mathcal{N}%
\models \left[ \exists \alpha \in \mathbf{Ord}^{\mathcal{N}}\ \exists
f:(V_{\alpha },\in )\overset{\cong }{\longrightarrow }(V_{\beta },\in )^{%
\mathcal{M}^{\ast }}\right] \right\} .
\end{equation*}%
$I$ and $J$ are nonempty since $0^{\mathcal{N}}\in I$ and $0^{\mathcal{M}%
^{\ast }}\in J$. $I$ and $J$ are clearly closed under predecessors. It is
also straightforward to show that $I$ and $J$ are closed under immediate
successors. To verify this, suppose $\alpha \in I$, then there is a
corresponding $\beta \in J$ and some $f_{\alpha }\in N$ such that:%
\begin{equation*}
\mathcal{N}\models \left[ f_{\alpha }:(V_{\alpha },\in )\overset{\cong }{%
\longrightarrow }(V_{\beta },\in )^{\mathcal{M}^{\ast }}\right] .
\end{equation*}%
We wish to show extend $f_{\alpha }$ to some $f_{\alpha +1}\in N$ such that: 
\begin{equation*}
\mathcal{N}\models \left[ f_{\alpha +1}:(V_{\alpha +1},\in )\overset{\cong }{%
\longrightarrow }(V_{\beta +1},\in )^{\mathcal{M}^{\ast }}\right] .
\end{equation*}%
Reasoning within $\mathcal{N}$, $f_{\alpha +1}(a)$ is defined for $a\in
V_{\alpha +1}$ as follows. Given $a\in V_{\alpha +1}$, $a\subseteq V_{\alpha
},$ so $\left\{ f_{\alpha }(x):x\in a\right\} $ is coded in $\mathcal{M}%
^{\ast }$ by a unique $b_{a}$ (using the assumptions $\mathcal{M}%
\trianglerighteq _{\mathrm{par}}\mathcal{N}$ and the existence of an $%
\mathcal{M}$-definable isomorphism between $\mathcal{M}$ and $\mathcal{M}%
^{\ast }),$ hence by the veracity of Replacement\footnote{%
In this step in the corresponding proof of solidity in the $\mathsf{KM}$%
-context, $\mathsf{CC}_{\mathsf{set}}$ is invoked within $(\mathcal{N},%
\mathcal{B})$ to be get hold of $f_{\alpha +1}\in \mathcal{B}.$} in $%
\mathcal{N}$ there is some $u\in N$ such that $\mathcal{N}\models u=\left\{
\left\langle a,b_{a}\right\rangle :a\in V_{\alpha +1}\right\} $, and thus
within $\mathcal{N}$ we can let $f_{\alpha +1}:=u.$ It is evident (again,
reasoning in $\mathcal{N},$ using $f_{\alpha }$, together with the
assumptions $\mathcal{M}\trianglerighteq _{\mathrm{par}}\mathcal{N}$ and the
existence of an $\mathcal{M}$-definable isomorphism between $\mathcal{M}$
and $\mathcal{M}^{\ast }$) that every $b\in V_{\beta +1}^{\mathcal{M}^{\ast
}}$ is in the image of $f_{\alpha +1}$. This shows that $\alpha +1\in I$ and 
$\beta +1\in J.$ So we have:\medskip 

\noindent (3) $\ \ I$ is an untopped\footnote{%
More explicitly: $I$ has no maximum element.} initial segment of $\mathbf{Ord%
}^{\mathcal{N}}$, and $J$ is an untopped initial segment of $\mathbf{Ord}^{%
\mathcal{M}^{\ast }}$.\medskip

\noindent By the rigidity of ordinals in $\mathsf{ZF}$, together with the
assumptions $\mathcal{M}\trianglerighteq _{\mathrm{par}}\mathcal{N}$ and the
existence of $\mathcal{M}$-definable isomorphism between $\mathcal{M}$ and $%
\mathcal{M}^{\ast }$, if $\alpha \in I$, then there is a unique $\beta \in J$%
, and a unique $f\in N$ such that $\mathcal{N}\models f:(V_{\alpha },\in )%
\overset{\cong }{\longrightarrow }(V_{\beta },\in )^{\mathcal{M}^{\ast }}$.
Therefore: $\medskip $

\noindent (4) $\ \ \mathcal{N}\models \forall \alpha \in I\ \exists !\beta
\in \mathbf{Ord}^{\mathcal{M}^{\ast }}$\ $\exists !f\ \in \mathbf{Ord}^{%
\mathcal{M}^{\ast }}\ f_{\alpha }:(V_{\alpha },\in )\overset{\cong }{%
\longrightarrow }(V_{\beta },\in )^{\mathcal{M}^{\ast }}.\medskip $

\noindent For $\alpha \in I,$ let $f_{\alpha }$ be the unique $f$ witnessing
(4). Note that (4) implies:$\medskip $

\noindent (5) $\ \ \mathcal{N}\models \forall \alpha ,\alpha ^{\prime }\in
I\ \left( \alpha \leq \alpha ^{\prime }\rightarrow f_{\alpha }\subseteq
f_{\alpha ^{\prime }}\right) .\medskip $

\noindent For $\alpha \in \mathbf{Ord}^{\mathcal{N}}$ and $\beta \in \mathbf{%
Ord}^{\mathcal{M}^{\ast }}$ let $\mathcal{N}_{\alpha }$ be the submodel of $%
\mathcal{N}$ whose universe is $N_{\alpha }:=\left\{ a\in N:\mathcal{N}%
\models a\in V_{\alpha }\right\} $ and $\mathcal{M}_{\beta }^{\ast }$ be the
submodel of $\mathcal{M}^{\ast }$ whose universe is $M_{\beta }^{\ast
}:=\left\{ b\in M^{\ast }:\mathcal{M}^{\ast }\models b\in V_{\beta }\right\} 
$. Next, let:$\medskip $

\begin{center}
$\mathcal{N}_{I}:=\bigcup\limits_{\alpha \in I}\mathcal{N}_{\alpha },$ and $%
\mathcal{M}_{J}^{\ast }:=\bigcup\limits_{\beta \in J}\mathcal{M}_{\beta
}^{\ast }.$
\end{center}

\noindent Define $F:\mathcal{N}_{I}\longrightarrow \mathcal{M}_{J}^{\ast }$,
given by: 
\begin{equation*}
F(a)=b\leftrightarrow \mathcal{N}\models \left[ \exists \alpha \in I\mathbf{%
\ }\left( a\in V_{\alpha }\wedge f_{\alpha }(a)=b\right) \right] .
\end{equation*}

\noindent $F$ is well-defined thanks to (4) and (5), so we have:$\medskip $

\noindent (6) \ \ $F:\mathcal{N}_{I}\overset{\cong }{\longrightarrow }%
\mathcal{M}_{J}^{\ast }$ and $F$ is $\mathcal{N}$-definable.$\medskip $

\noindent At this point, we distinguish four cases:$\medskip $

\noindent Case 1. $I$ is bounded in $\mathbf{Ord}^{\mathcal{N}}$, and $J$ is
bounded in $\mathbf{Ord}^{\mathcal{M}^{\ast }}.\medskip $

\noindent Case 2. $I$ is bounded in $\mathbf{Ord}^{\mathcal{N}}$, and $J$ $=%
\mathbf{Ord}^{\mathcal{M}^{\ast }}.\medskip $

\noindent Case 3. $I=\mathbf{Ord}^{\mathcal{N}}$, and $J$ is bounded in $%
\mathbf{Ord}^{\mathcal{M}^{\ast }}.\medskip $

\noindent Case 4. $I=\mathbf{Ord}^{\mathcal{N}}$, and $J=\mathbf{Ord}^{%
\mathcal{M}^{\ast }}.\medskip $

\noindent If Case 4 holds, then by (6) $G$ is an $\mathcal{N}$-definable
isomorphism between $\mathcal{N}$ and $\mathcal{M}^{\ast }$. Thus the proof
of the theorem will be complete once we rule out Cases 1 through 3.$\medskip 
$

\noindent If Case 1 holds then let $\overline{\alpha }:=\min (\mathbf{Ord}^{%
\mathcal{N}}\backslash I\mathbf{),}$ and $\overline{\beta }:=\min (\mathbf{%
Ord}^{\mathcal{M}^{\ast }}\backslash J\mathbf{)}$. The well-definedness of $%
\overline{\alpha }$ is obvious, and the well-definedness of $\overline{\beta 
}$ is assured by (2). By (3), $\overline{\alpha }$ is a limit ordinal of $%
\mathcal{N}$ and $\overline{\beta }$ is a limit ordinal of $\mathcal{M}%
^{\ast }$. By the veracity of Replacement\footnote{%
In this step in the corresponding proof of solidity in the $\mathsf{KM}$%
-context, $\mathsf{CC}_{\mathsf{set}}$ is invoked within $(\mathcal{N},%
\mathcal{B})$ to be get hold of $f_{\overline{\alpha }}\in \mathcal{B}.$} in 
$\mathcal{N}$, there is some $v\in N$ such that $\mathcal{N}\models \left[
v=\left\{ f_{\alpha }:\alpha \in I\right\} \right] $, and therefore there is 
$f_{\overline{\alpha }}\in N$ such that $\mathcal{N}\models \left[ f_{%
\overline{\alpha }}=\bigcup\limits_{\alpha \in \overline{\alpha }}f_{\alpha }%
\right] .$ Since $\overline{\alpha }$ and $\overline{\beta }$ are limit
ordinals in their respective models, this show that $f_{\overline{\alpha }}$
is an isomorphism between $\mathcal{N}_{\overline{\alpha }}$ and $\mathcal{M}%
_{\overline{\beta }}^{\ast }$. Thus $\overline{\alpha }\in I$, which
contradicts the definition of $\overline{\alpha }.$ \medskip

\noindent Now suppose Case 2 holds, then by (6) $F$ is a $\mathcal{N}$%
-definable isomorphism between $\mathcal{N}_{\overline{\alpha }}$ and $%
\mathcal{M}^{\ast }$, where $\overline{\alpha }$ is as in Case 1(b). Since
the satisfaction predicate for $\mathcal{N}_{\overline{\alpha }}$ is $%
\mathcal{N}$-definable, we can use the assumptions $\mathcal{M}%
\trianglerighteq _{\mathrm{par}}\mathcal{N}$, and $\mathcal{M}\models \left[
i_{0}^{-1}:\mathcal{M}^{\ast }\overset{\cong }{\longrightarrow }(\mathbf{V}%
,\in )\right] $ to show that the satisfaction predicate for $\mathcal{M}$ is 
$\mathcal{M}$-definable, which contradicts Tarski's Undefinability of Truth
Theorem.$\medskip $

\noindent Finally, if Case 3 holds then by (7) $F$ is a $\mathcal{N}$%
-definable isomorphism between $\mathcal{N}$ and $\mathcal{M}_{\overline{%
\beta }}^{\ast }$, where $\overline{\beta }$ is as in Case 1(b). Using the
assumptions $\mathcal{M}\trianglerighteq _{\mathrm{par}}\mathcal{N}$, and $%
\mathcal{M}\models \left[ i_{0}^{-1}:\mathcal{M}^{\ast }\overset{\cong }{%
\longrightarrow }(\mathbf{V},\in )\right] $, it is evident that $N$ is a set
in $\mathcal{M}$, which contradicts (1). \hfill $\square $\medskip

We are now ready to discuss the gap in the proof of solidity of $\mathsf{KM}$
presented in \cite[Theorem 2.7]{Enayat-Visserian}. Recall that models of $%
\mathsf{KM}$ can be represented as two-sorted structures of the form $\left( 
\mathcal{M},\mathcal{A}\right) $, where $\mathcal{M}\models \mathsf{ZF}$; $%
\mathcal{A}$ is a collection of subsets of $M$; and $\left( \mathcal{M},%
\mathcal{A}\right) $ satisfies the full comprehension scheme. \medskip

\noindent \textbf{2.2.~Definition.} Let $\mathcal{L}_{\mathsf{KM}}$ be the
two-sorted language of $\mathsf{KM}$.\medskip

\noindent \textbf{(a)} $\mathsf{CC}$ (Class Collection) consists of the
universal closure of formulae of the form: 
\begin{equation*}
\left[ \forall x\ \exists X\ \varphi (x,X)\right] \rightarrow \left[ \exists
Y\ \forall x\ \varphi (x,(Y)_{x})\right] ,
\end{equation*}%
where $\varphi (x,X)$ is an $\mathcal{L}_{\mathsf{KM}}$-formula in which $Y$
does not occur free, and is allowed to have set or class parameters. Here $%
(Y)_{x}$ is the `$x$-th cross section' of the subclass of $\mathbf{V}^{2}$
coded by $Y$, i.e., 
\begin{equation*}
(Y)_{x}=\{y\in \mathbf{V}:\left\langle x,y\right\rangle \in Y\},
\end{equation*}%
where $\left\langle x,y\right\rangle $ is the Kuratowski ordered pair of $x$
and $y$.\footnote{%
We have followed \cite{GHJ} in adopting $\mathsf{CC}$ to refer to this
scheme in lieu of $\mathsf{AC}$\textit{\ }(used, e.g., in \cite%
{Antos-Friedman} and \cite{Enayat-Lelyk-JPM}), which in a context involving
sets, runs the risk of getting confused with the axiom of choice.}\medskip

\noindent \textbf{(b)} $\mathsf{CC}_{\mathsf{set}}$ is the fragment of $%
\mathsf{CC}$ that consists of the universal closure of formulae of the form: 
\begin{equation*}
\left[ \forall x\in s\ \exists X\ \varphi (x,X)\right] \rightarrow \left[
\exists Y\ \forall x\in s\ \varphi (x,(Y)_{x})\right] .
\end{equation*}

\noindent \textbf{(c)} $\mathsf{CR}$ (Class Replacement) is the fragment of $%
\mathsf{CC}$ that consists of the universal closure of formulae of the form: 
\begin{equation*}
\left[ \forall x\ \exists !X\ \varphi (x,X)\right] \rightarrow \left[
\exists Y\ \forall x\ \varphi (x,(Y)_{x})\right] .
\end{equation*}

\noindent \textbf{2.3.~Remark.} As shown by Gitman, Hamkins, and Karagila 
\cite{Gitam-et-al(FODOR)}, $\mathsf{CC}$ is not provable in $\mathsf{KM}$.%
\footnote{%
The analogue of $\mathsf{CC}$ is the context of $\mathsf{Z}_{2}$ (second
order arithmetic) is commonly referred to $\mathsf{AC}$. It has long been
known that $\mathsf{Z}_{2}$ is capable of proving the fragment $\Sigma
_{2}^{1}$-$\mathsf{AC}$ of $\mathsf{AC}$ (using an appropriate
implementation of the Shoenfield Absoluteness in $\mathsf{Z}_{2}$; see
Theorem VII.6.9.1 of Simpson's monograph \cite{Steve Book}), but not full $%
\mathsf{AC}$. Indeed the standard model of $\mathsf{Z}_{2}$ of the classical
Feferman-Levy model of $\mathsf{ZF}$ in which $\aleph _{1}$ is a countable
union of countable sets has the property that $\mathsf{AC}$ fails for some $%
\Pi _{2}^{1}$-formula $\varphi $; see \cite[Theorem 8]{Levy}. In contrast, a
straightforward induction on the variable $s$ shows that $\mathsf{Z}_{2}$ is
capable of proving the arithmetical analogue of $\mathsf{CC}$ whose
instances are of the form:
\par
\begin{center}
$\left[ \forall x<y\ \exists X\ \varphi (x,X)\right] \rightarrow \left[
\exists Y\ \forall x<y\ \varphi (x,(Y)_{x})\right] .$%
\end{center}
} This fact was subsequently refined in the work of Gitman, Hamkins, and
Johnstone \cite{GHJ}, where it is shown that $\mathsf{KM}$ is incapable of
even proving those instances of $\mathsf{CC}_{\mathsf{set}}$ in which $%
s=\omega $; see Theorem 18, 19, and 20 of \cite{GHJ}. In contrast, $\mathsf{%
CR}$ is provable in $\mathsf{KM}$ since if $\forall x\ \exists !X\ \varphi
(x,X)$ holds, then we can use comprehension applied to $\psi (x,y):=\exists
X\ (\varphi (x,X)\wedge y\in X)$ to obtain the class $Y:=\left\{
\left\langle x,y\right\rangle :\psi (x,y)\right\} .$ Then, by design, $%
\forall x\ \varphi (x,(Y)_{x})$ holds, as desired. \footnote{%
I am grateful to Vika Gitman for pointing out the fact that $\mathsf{CR}$ is
provable in $\mathsf{KM.}$}\medskip

\noindent \textbf{2.4.~Theorem.} $\mathsf{KM+CC}$ \textit{does not eliminate
imaginaries}.\footnote{%
A theory $T$ is said to eliminate imaginaries, if for each definable
equivalence relation $E$, then there is a definable function $f$ such that
for all $x$ and $y$, $xEy$ $\leftrightarrow f(x)=f(y).$ In $\mathsf{ZF}$
such a function can be readily described by the so-called `Scott trick',
which takes advantage of the stratification of the universe into the
well-ordered family of sets of the form $V_{\alpha }$, as $\alpha $ ranges
over the ordinals. For more about this concept, see Section 4.4 of Hodges'
textbook \cite{Hodges Book}.}\medskip 

\noindent \textbf{Proof.} Let $\mathsf{ZF}^{-}$ be the result of eliminating
the powerset axiom from the formulation of $\mathsf{ZF}$ in which the
separation scheme and the collection scheme take the place of the
replacement scheme. Also, let:\medskip

\begin{center}
$T_{1}:=\mathsf{KM+CC}$, and $T_{2}:=\mathsf{ZF}^{-}+\exists \kappa $ ($%
\kappa $ is inaccessible and $\forall x\ \left\vert x\right\vert \leq \kappa 
$)$.$ \medskip
\end{center}

\noindent The result follows in light of the following facts:

\begin{enumerate}
\item[$(1)$] $T_{1}$ and $T_{2}$ are sequential theories. This is trivial
once one knows the definition of sequentiality. Roughly speaking, a theory
is sequential if it supports a modicum of coding machinery to have access to
a `$\beta $-function' that codes finite sequences of objects in the domain
of discourse.\footnote{%
More precisely, a theory $T$ is sequential if there is a formula $N(x)$,
together with appropriate formulae providing interpretations of equality,
and the operations of successor, addition, and multiplication for elements
satisfying $N(x)$ such that $T$ proves the translations of the axioms of $%
\mathsf{Q}$ (Robinson's arithmetic) when relativized to $N(x)$; and
additionally, there is a formula $\beta (x,i,w)$ (whose intended meaning is
that $x$ is the $i$-th element of a sequence $w$) such that $T$ proves that
every sequence can be extended by any given element of the domain of
discourse. }

\item[$(2)$] $T_{1}$ is bi-interpretable with $T_{2}$, where the witnessing
interpretation of $T_{1}$ in $T_{2}$ is identity preserving, but the
witnessing interpretations of $T_{2}$ in $T_{1}$ is not identity preserving.
This is a well-known fact, the main ingredients of whose proof were
developed by Mostowski and Marek. For a perspicuous exposition, see Chapter
2 of William's doctoral dissertation.

\item[$(3)$] If two sequential theories are bi-interpretable via a pair of
equality-preserving interpretations, then they are definitionally
equivalent. This result is due to Friedman and Visser \cite{Harvey+Albert}.

\item[$(4)$] $T_{1}$ and $T_{2}$ are not definitionally equivalent; this was
established in \cite[Theorem 16]{Enayat-Lelyk-JPM}.\footnote{%
Chen and Meadows \cite{Chen+Meadows} gave a simpler proof than the one
offered in \cite{Enayat-Lelyk-JPM} for the failure of definitional
equivalence of the arithmetical counterparts of $T_{1}$ and $T_{2}.$ Their
proof strategy can be used to provide a simpler proof of the failure of
definitional equivalence of $T_{1}$ and $T_{2}$ as well.}
\end{enumerate}

\noindent More explicitly, if $\mathsf{KM+CC}$ could eliminate imaginaries,
then by (1), (2), and (3), $T_{1}$ and $T_{2}$ would have to be
definitionally equivalent, which contradicts (4). \hfill $\square $\medskip

\noindent \textbf{2.5.~Remark. }In light of Theorem 2.4, in treating
interpretations in the $\mathsf{KM}$ context, have to allow for equality to
be interpreted as an equivalence relation. For this purpose, let us clarify
what is meant by a definable isomorphism between two interpreted structures.
Suppose $\mathcal{M}$ is structure that interprets structures $\mathcal{N}%
_{i}$ for $i\in \{1,2\}$, using the domain formulae $\delta _{i}$, and
definable equivalence relations $E_{i}$. An $\mathcal{M}$-definable
isomorphism between $\mathcal{N}_{0}$ and $\mathcal{N}_{1}$ is an $\mathcal{M%
}$-definable relation $R\subseteq \delta _{0}^{\mathcal{M}}\times \delta
_{1}^{\mathcal{M}}$ that satisfies the following properties:

\begin{enumerate}
\item The domain and the codomain of $R$ are $\delta _{0}^{\mathcal{M}}$ and 
$\delta _{1}^{\mathcal{M}}$, respectively.

\item If $x_{0}Rx_{1}$ and $y_{0}E_{0}x_{0}$ and $y_{1}E_{1}x_{1}$, then $%
y_{0}Ry_{1}$.

\item The function $F:N_{0}\rightarrow N_{1}$, defined by $%
F([x_{0}]_{E_{0}})=[x_{1}]_{E_{1}}\Longleftrightarrow x_{0}Rx_{1}$, is an
isomorphism between $\mathcal{N}_{0}$ and $\mathcal{N}_{1}.$\medskip
\end{enumerate}

\begin{itemize}
\item If $F$ and $R$ are as above, then we will refer to the relation $R$ as
a\textit{\ representation} of the function $F$.
\end{itemize}

\noindent \textbf{2.6.~Remark. }In our two-sorted set-up for $\mathsf{KM}$,
if $\left( \mathcal{M},\mathcal{A}\right) $ and $\left( \mathcal{N},\mathcal{%
B}\right) $ are models of $\mathsf{KM}$, when asserting that 
\begin{equation*}
\left( \mathcal{M},\mathcal{A}\right) \trianglerighteq _{\mathrm{par}}\left( 
\mathcal{N},\mathcal{B}\right) ,
\end{equation*}%
we allow the entanglement of the sorts in the interpretation. Coupled with
Theorem 2.4, all we can infer about $\mathcal{N}$ on the basis of the above
displayed assertion is that there are: $(i)$ an $\left( \mathcal{M},\mathcal{%
A}\right) $-definable subset $D\subseteq M\cup \mathcal{A}$, $(ii)$ some $%
\left( \mathcal{M},\mathcal{A}\right) $-definable equivalence relation $E$
on $D$, and $(iii)$ some $\left( \mathcal{M},\mathcal{A}\right) $-definable
binary relation $F$ on $D$ such that:

\begin{equation*}
\mathcal{N}\cong \left( D/E,F\right) .
\end{equation*}%
As we shall explain, without loss of generality, in $(i)$ we can assume that 
$D\subseteq \mathcal{A}$, essentially because of the fact that there is a
parameter-free $\left( \mathcal{M},\mathcal{A}\right) $-definable bijection
between $\mathcal{A}$ and $M\cup \mathcal{A}$, we can assume without loss of
generality that $D\subseteq \mathcal{A}.$\footnote{%
This is trivially true in the one-sorted set-up for \textsf{KM}, where all
sets are counted as classes. In contrast, in the two-sorted formulation,
there is no overlap between objects of different sorts. Using the bijection $%
g$ described in Remark 2.2, one can readily check that the two formulations
of $\mathsf{KM}$ are bi-interpretable; the argument can be readily carried
out for the much weaker theory $\mathsf{GB}$ (G\"{o}del-Bernays theory of
classes) as well.} To see this, we first observe that within $\mathsf{ZF}$
the universe $\mathbf{V}$ can be written as $\mathbf{V}_{1}\cup \mathbf{V}%
_{2}$, where $\mathbf{V}_{1}$ and $\mathbf{V}_{2}$ are disjoint and there is
a definable bijection between $\mathbf{V}$ and $\mathbf{V}_{i}$ for $i\in
\{1,2\}.$ This can be accomplished, e.g., by letting $\mathbf{V}_{1}$ to be
the set of all singletons, and $\mathbf{V}_{2}:=\mathbf{V}\backslash \mathbf{%
V}_{1}.$ Clearly the function $f_{1}:\mathbf{V}\rightarrow \mathbf{V}_{1}$
given by $f_{1}(x)=\{x\}$ is a definable bijection between $f_{1}:\mathbf{V}%
\rightarrow \mathbf{V}_{1}.$ Using the class-version of the Schr\"{o}%
der-Bernstein theorem, there is also a bijection $f_{2}:\mathbf{V}%
\rightarrow \mathbf{V}_{2},$ since the inclusion map is an injection of $%
\mathbf{V}_{2}$ into $\mathbf{V}$, and the map $x\longmapsto \{x,\{x\}\}$ is
an injection of $\mathbf{V}$ into $\mathbf{V}_{2}.$ One can then take
advantage of $f_{1}$ and $f_{2}$ together with the injection $g$ from $%
\mathbf{V}$ to the collection of classes $\mathbf{C}$ given by $g(x)=X$,
where $X$ only contains $x,$ to construct a definable bijection $g:\mathbf{%
C\rightarrow V}\cup \mathbf{C}$ within $\mathsf{KM}$ between $\mathbf{V}\cup 
\mathbf{C}$ and $\mathbf{C}$ by:%
\begin{equation*}
g(x)=\left\{ 
\begin{tabular}{ll}
$x$ & \textrm{if}$\ x\notin i(\mathbf{V})$ \\ 
$f_{1}(i(x))$ & \textrm{if}$\ x\in i(\mathbf{V}_{1})$ \\ 
$f_{2}(i(x))$ & \textrm{if}$\ x\in i(\mathbf{V}_{2})$%
\end{tabular}%
\right. .
\end{equation*}

\noindent \textbf{2.7.~Theorem. }$\mathsf{KM+CC}_{\mathsf{set}}$ \textit{is
solid.}\medskip

\noindent \textbf{Proof. }Suppose $\left( \mathcal{M},\mathcal{A}\right) $, $%
\left( \mathcal{M}^{\ast },\mathcal{A}^{\ast }\right) $, and $\left( 
\mathcal{N},\mathcal{B}\right) $ are models of $\mathsf{KM}$ such that:%
\begin{equation*}
\left( \mathcal{M},\mathcal{A}\right) \trianglerighteq _{\mathrm{par}}\left( 
\mathcal{N},\mathcal{B}\right) \trianglerighteq _{\mathrm{par}}\left( 
\mathcal{M}^{\ast },\mathcal{A}^{\ast }\right) ,
\end{equation*}

\noindent and there is an $\left( \mathcal{M},\mathcal{A}\right) $-definable
isomorphism\textit{\ }$\widehat{i_{0}}:\left( \mathcal{M},\mathcal{A}\right) 
\overset{\cong }{\rightarrow }\left( \mathcal{M}^{\ast },\mathcal{A}^{\ast
}\right) .$ As noted by Vika Gitman, if the interpretations at work are
identity preserving, then the proof of \cite[Theorem 2.7]{Enayat-Visserian}
is fine and does not need the extra $\mathsf{CC}_{\mathsf{set}}$ assumption,
basically because of the provability of Class Replacement in $\mathsf{KM}$
(as explained in Remark 2.3). Indeed, it is sufficient for the interpetation
witnessing $\left( \mathcal{N},\mathcal{B}\right) \trianglerighteq _{\mathrm{%
par}}\left( \mathcal{M}^{\ast },\mathcal{A}^{\ast }\right) $ to be identity
preserving.\medskip 

\noindent To handle the general case where the interpretations at work do
not necessarily translate equality as equality, suppose, furthermore, that $%
\mathsf{CC}_{\mathsf{set}}$ holds in $\left( \mathcal{N},\mathcal{B}\right) $%
. We wish to repeat the argument of the $\mathsf{ZF}$ case to show that
there is an $\left( \mathcal{N},\mathcal{B}\right) $-definable isomorphism $%
j:\mathcal{N}\rightarrow \mathcal{M}^{\ast }.$ Note that if we succeed in
doing so, then using the argument in the proof of \cite[Theorem 2.7]%
{Enayat-Visserian} $j$ can be naturally lifted to an $\left( \mathcal{N},%
\mathcal{B}\right) $-definable isomorphism 
\begin{equation*}
\widehat{j}:\left( \mathcal{N},\mathcal{B}\right) \overset{\cong }{%
\rightarrow }\left( \mathcal{M}^{\ast },\mathcal{A}^{\ast }\right) .
\end{equation*}

\noindent In the $\mathsf{KM}$ context, the counterpart of statement (1) of
the proof of Theorem 2.1 asserts that the sentence $\lnot \exists x\forall
y(y\in x\leftrightarrow y\in N)$ holds in $\mathcal{(M},\mathcal{A)}$, and
the counterpart of statement (2) asserts that $\in ^{\mathcal{M}^{\ast }}$
is well-founded as viewed by $\left( \mathcal{N},\mathcal{B}\right) $ in the
strong sense that any nonempty subset of $M^{\ast }$ that is definable in $%
\left( \mathcal{N},\mathcal{B}\right) $ has an $\in ^{\mathcal{M}^{\ast }}$%
-least element. Both of the aforementioned counterparts are readily provable
in $\mathsf{KM}$ (without the use of $\mathsf{CC}_{\mathsf{set}}$) using
essentially the same proofs as in the $\mathsf{ZF}$ case. Then, we define
the relevant initial segments $I$ of $\mathbf{Ord}^{\mathcal{N}}$ and $J$ of 
$\mathbf{Ord}^{\mathcal{M}^{\ast }}$ as follows:%
\begin{equation*}
I=\left\{ \alpha \in \mathbf{Ord}^{\mathcal{N}}\mid \left( \mathcal{N},%
\mathcal{B}\right) \models \left[ \exists \beta \in \mathbf{Ord}^{\mathcal{M}%
^{\ast }}\ \exists f\in \mathcal{B\ \ }f:(V_{\alpha },\in )\overset{\cong }{%
\longrightarrow }(V_{\beta },\in )^{\mathcal{M}^{\ast }}\right] \right\} ,
\end{equation*}%
and 
\begin{equation*}
J=\left\{ \beta \in \mathbf{Ord}^{\mathcal{M}^{\ast }}\mid \left( \mathcal{N}%
,\mathcal{B}\right) \models \left[ \exists \alpha \in \mathbf{Ord}^{\mathcal{%
N}}\ \ \exists f\in \mathcal{B\ \ }f:(V_{\alpha },\in )\overset{\cong }{%
\longrightarrow }(V_{\alpha },\in )^{\mathcal{M}^{\ast }}\right] \right\} .
\end{equation*}

\noindent Note that in contrast with the proof of Theorem 2.1, in the
definition of $I$ and $J$ we have to allow $f$ to be a \textit{class} in the
sense of $\left( \mathcal{N},\mathcal{B}\right) $, but $f$ need not
literally be the graph of a \textit{function} in the sense of $\left( 
\mathcal{N},\mathcal{B}\right) ,$ since equality of $\mathcal{M}^{\ast }$
might be translated as an $\left( \mathcal{N},\mathcal{B}\right) $-definable
equivalence relation $E$. Instead, $f$ should be function-like if\textit{\ }$%
E$\textit{\ }serves as equality for outputs, i.e., $\left( \mathcal{N},%
\mathcal{B}\right) $ should satisfy: 
\begin{equation*}
\forall x\forall y\ \left[ \left( (x,y)\in f\wedge (x,z)\in f\right)
\rightarrow yEz\right] .
\end{equation*}%
Examining the proof strategy of Theorem 2.1, we can readily carry out all
the steps of the proof of Theorem 2.1 within $\mathsf{KM}$, with the
provisos\ below. We assume that $E$ is the $\left( \mathcal{N},\mathcal{B}%
\right) $-definable equivalence relation that is the interpretation of
equality of $\mathcal{M}^{\ast }$.

\begin{itemize}
\item In the proof of the analogue of statement (3) of the proof of Theorem
2.1, to show that $f_{\alpha +1}\in \mathcal{B}$, the veracity of $\mathsf{CC%
}_{\mathsf{set}}$ in $\left( \mathcal{N},\mathcal{B}\right) $ plays the role
of $\mathcal{N}$ satisfying Replacement.

\item Statement (4) of the proof of Theorem 2.1 asserting that $\beta $ and $%
f$ are unique, translates in $\left( \mathcal{N},\mathcal{B}\right) $ as: If 
$\alpha \in I,$ and $\beta \in J$, and $f\in \mathcal{B}$ such that $\left( 
\mathcal{N},\mathcal{B}\right) \models \left[ \mathcal{\ }f:(V_{\alpha },\in
)\overset{\cong }{\longrightarrow }(V_{\beta },\in )^{\mathcal{M}^{\ast }}%
\right] ,$ and $\beta _{1}$ and $\beta _{2}$ are elements of $\mathcal{B}$
that serve as representatives of $\beta $, and $f_{1}$ and $f_{2}$ are
elements of $\mathcal{B}$ that serve as representatives of $f$, then: 
\begin{equation*}
\left( \mathcal{N},\mathcal{B}\right) \models \left[ \beta _{1}E\beta
_{2}\wedge \forall x_{1}\forall x_{2}\forall y_{1}\forall y_{2}\left(
x_{1}Ex_{2}\rightarrow \left( \left( \left\langle x_{1},y_{1}\right\rangle
\in f_{1}\wedge \left\langle x_{2},y_{2}\right\rangle \in f_{2}\right)
\rightarrow y_{1}Ey_{2}\right) \right) \right] .
\end{equation*}

\item Similarly, statement (6) of the proof of Theorem 2.1 that asserts that 
$f_{\alpha }\subseteq f_{\alpha ^{\prime }}$ for $\alpha ,\alpha ^{\prime
}\in I$ with $\alpha \leq \alpha ^{\prime }$, is translated within $\left( 
\mathcal{N},\mathcal{B}\right) $ as: $\left[ (x,y)\in f_{\alpha }\right]
\wedge \left[ (x,z)\in f_{\beta }\right] \rightarrow yEz.$

\item $\mathcal{N}_{I}$ and $\mathcal{M}_{J}^{\ast }$ are defined as in the
proof of Theorem 2.1, and the $(\mathcal{N},\mathcal{B}$)-definable
isomorphism $F:\mathcal{N}_{I}\longrightarrow \mathcal{M}_{J}^{\ast }$ is
given by the following, where $\left\langle \left\langle a,b\right\rangle
\right\rangle $ is the ordered pair appropriate for this context.\footnote{%
Here, $a$ is a set in $(\mathcal{N},\mathcal{B})$, but $b$ is allowed to be
in $\mathcal{B}$, so we cannot use the Kuratowski ordered pair formation $%
\left\langle a,b\right\rangle $, instead we can use $\left\langle
\left\langle a,b\right\rangle \right\rangle :=$ $\{\left\langle
x,0\right\rangle :x\in a\}\cup \{\left\langle y,1\right\rangle :y\in b\}.$} 
\begin{equation*}
F:=\left\{ \left\langle \left\langle a,b\right\rangle \right\rangle \mid
\exists f\in \mathcal{F\ }\left( \mathcal{N},\mathcal{B}\right) \models
\left\langle \left\langle a,b\right\rangle \right\rangle \in f\right\} ,
\end{equation*}%
where:%
\begin{equation*}
\mathcal{F}:=\left\{ f\in \mathcal{B}\mid \left( \mathcal{N},\mathcal{B}%
\right) \models \exists \alpha \in I\ \exists \beta \in J\mathcal{\ \ }%
f:(V_{\alpha },\in )\overset{\cong }{\longrightarrow }(V_{\alpha },\in )^{%
\mathcal{M}^{\ast }}\right\} .\footnote{%
In Cases 3 and 4, when $I=\mathbf{Ord}^{\mathcal{N}}$, there is no
requirement that there is some $Z\in \mathcal{B}$ that codes $F$, since we
only require that $F$ is $\left( \mathcal{N},\mathcal{B}\right) $-definable.
However, such a $Z$ exists if the full $\mathsf{CC}$ holds in $\left( 
\mathcal{N},\mathcal{B}\right) $.}
\end{equation*}

\item In the proof of the impossibility of Case 1, to show that there is a
representation of $f_{\overline{\alpha }}$ in $\mathcal{B}$, the assumption
that $\mathsf{CC}_{\mathsf{set}}$ holds in $\left( \mathcal{N},\mathcal{B}%
\right) $ plays the role of Replacement holding in $\mathcal{N}$.
\end{itemize}

\hfill $\square $ \medskip

\noindent \textbf{2.8.~Remark.~}Recall from \cite{Enayat-Visserian} that $%
\mathsf{KM}_{n}$ is the $n$-th order Kelley-Morse theory of classes; in this
notation $\mathsf{KM}$ is $\mathsf{KM}_{1}.$ A reasoning similar to the
proof of Theorem 2.7 shows that $\mathsf{KM}_{n}+\mathsf{CC}_{n\mathsf{,set}%
} $ is solid; where $\mathsf{CC}_{n\mathsf{,set}}$ is the analogue of $%
\mathsf{CC}_{\mathsf{set}}$ for $\mathsf{KM}_{n}$, which consists of the
universal closure of formulae of the form 
\begin{equation*}
\left[ \forall x\in s\ \exists ^{n+1}X\ \varphi (x,X)\right] \rightarrow %
\left[ \exists ^{n+1}Y\ \forall x\in s\ \varphi (x,(Y)_{x})\right] ,
\end{equation*}%
where $x$ and $s$ are first order object, $X$ and $Y$ are $\left( n+1\right) 
$-th order objects, and $\varphi (x,X)$ is a formula (with suppressed
parameters) in the language of $\mathsf{KM}_{n}.$


\medskip

\begin{thebibliography}{E\L -1}
\bibitem[AF]{Antos-Friedman} C.~Antos and S.D.~Friedman, \textit{Hyperclass
forcinng in Morse-Kelley class theor}y, \textbf{Journal of Symbolic Logic},
vol.~82 (2017), pp.~549-575.

\bibitem[CM]{Chen+Meadows} J.~Chen and T.~Meadows, \textit{Teasing apart
definitional equivalence, }arXiv:2508.03956.

\bibitem[E]{Enayat-Visserian} A.~Enayat, \textit{Variations on a Visserian
theme}, in \textbf{Liber Amicorum Alberti. A Tribute to Albert Visser},
edited by J.~van Eijk, R.~Iemhoff, and J.~Joosten, College Publications,
pp.~99-110, 2016; {\small https://doi.org/10.48550/arXiv.1702.07093}.

\bibitem[E\L -1]{Enayat-Lelyk-JPM} A. Enayat and M. \L e\l yk, \textit{%
Categoricity-like properties in the first order realm}, \textbf{Journal for
the Philosophy of Mathematics}, vol.~pp.~1:63--98, (2024); {\small %
https://doi.org/10.36253/jpm-2934}

\bibitem[E\L -2]{Enayat-Lelyk-correction} A. Enayat and M. \L e\l yk, 
\textit{Corrigendum \& Addendum to: Categoricity-like properties in the
first order realm}, to appear.

\bibitem[FV]{Harvey+Albert} H.~Friedman and A.~Visser, \textit{When
Bi-Interpretability Implies Synonymy}, \textbf{Review of Symbolic Logic}
vol.~18, pp.~971-990 (2025).

\bibitem[GHJ]{GHJ} V.~Gitman, J.D.~Hamkins, and T.~Johnstone, \textit{Class
Choice and the surprising weakness of Kelley-Morse set theory},
arXiv:2601.23165.

\bibitem[GHK]{Gitam-et-al(FODOR)} V.~Gitman, J.D.~Hamkins, and A.~Karagila, 
\textit{Kelley-Morse set theory does not prove the class Fodor principle}, 
\textbf{Fund.}~\textbf{Math.}~vol.~254 (2021), pp.~133--154.

\bibitem[H]{Hodges Book} W. Hodges, \textbf{Model Theory}, Cambridge
University Press, 1993.

\bibitem[L]{Levy} A.~Levy, \textit{Definability in axiomatic set theory. II}%
, in \textbf{Mathematical Logic and Foundations of Set Theory}
(Proc.~Internat.~Colloq., Jerusalem, 1968), North-Holland, Amsterdam, 1970,
pp.~129--145.

\bibitem[M]{Meadows} T.~Meadows, \textit{Foundation with Imagination},
arXiv:2601.20057.

\bibitem[S]{Steve Book} S.~Simpson, \textbf{Subsystems of Second Order
Arithmetic}, Springer, Heidelberg 1999.

\bibitem[W]{Kameryn} K.~J.~Williams, \textbf{The Structure of Models of
Second-order Set Theories}, Doctoral Dissertation, 2018, arXiv:1804.09526.
\end{thebibliography}
\end{document}